\documentclass[12pt,reqno]{amsart}

\usepackage{undertilde, slashbox}
\usepackage[final
]{showkeys}

\usepackage{AKstyle}

\newcommand{\vertiii}[1]{{\left\vert\kern-0.25ex\left\vert\kern-0.25ex\left\vert #1
    \right\vert\kern-0.25ex\right\vert\kern-0.25ex\right\vert}}

\numberwithin{equation}{section}

\newcommand{\mult}{\text{mult}}

\usepackage{caption}
\usepackage[labelfont=rm]{subcaption}

\usepackage[pagebackref=true, colorlinks]{hyperref}

\hypersetup{pdffitwindow=true,linkcolor=blue,citecolor=blue,urlcolor=blue,menucolor=blue}

\usepackage{comment}

\begin{document}

\author{Jean Bourgain}
\thanks{Bourgain is partially supported by NSF grant DMS-1301619.}
\email{bourgain@math.ias.edu}
\address{IAS, Princeton, NJ}
\author{Alex Kontorovich}
\thanks{Kontorovich is partially supported by
an NSF CAREER grant DMS-1254788 and  DMS-1455705, an NSF FRG grant DMS-1463940, an Alfred P. Sloan Research Fellowship, and a BSF grant.}
\email{alex.kontorovich@rutgers.edu}
\address{Rutgers University, New Brunswick, NJ}
\author{Michael Magee}
\thanks{Magee is partially supported by NSF grant DMS-1128155.}
\email{mmagee@math.ias.edu}
\address{IAS, Princeton, NJ}

\title{Thermodynamic expansion to arbitrary moduli}

\begin{abstract}
We extend the thermodynamic expansion results in \cite{BourgainGamburdSarnak2011,  MageeOhWinter2015} from square-free to arbitrary moduli by developing a novel decoupling technique and applying   \cite{BourgainVarju2012}.
\end{abstract}
\date{\today}
\maketitle
\tableofcontents

\section{Statements}

In this short note, we use the ``modular'' expansion of \cite{BourgainVarju2012}, valid for arbitrary moduli, to extend the ``archimedean''-thermodynamic expansion results in \cite{BourgainGamburdSarnak2011, MageeOhWinter2015} from square-free to arbitrary moduli.
\begin{thm}\label{thm:1}
Let $\G$ be a finitely-generated, Zariski dense, Schottky (that is, free, convex-cocompact) subgroup of $\SL_{2}(\Z)$, and let $\gd\in(0,1)$ be its critical exponent.
For an integer $q$, let $\G(q):=\{\g\in\G:\g\equiv I(\mod q)\}$.
 Then there is an $\vep>0$ and $q_{0}\ge1$ such that, for all integers $q$ coprime to $q_{0}$, the resolvent of the Laplace operator
 $$
 R_{\G(q)}=(\gD-s(1-s))^{-1}:C_{c}^{\infty}(\G(q)\bk\bH)\to C^{\infty}(\G(q)\bk\bH)
 $$ 
 is holomorphic in the strip
 $$
 \Re(s)>\gd-\vep,
 $$
 except for a simple pole at $s=\gd$.
\end{thm}
This extends the statements of  \cite[Theorems 1.4 and 1.5]{BourgainGamburdSarnak2011} and \cite[Theorem 1.3]{OhWinter2014} to arbitrary moduli $q$; see also the discussion below \cite[Theorem 1.1]{MageeOhWinter2015}.

In a similar way, we deal with semigroups.
\begin{thm}\label{thm:2}
The statement of \cite[Corollary 1.2]{MageeOhWinter2015} holds
 when specialized to the ``Zaremba'' (or continued fraction) setting of \cite[\S 6.1]{MageeOhWinter2015}, without the restriction that the modulus  $q$ be square-free.
\end{thm}

In particular, this justifies Remark 8.8 in \cite{BourgainKontorovich2014}. We expect analogous arguments will also prove the uniform exponential mixing result in
\cite[Theorem 1.1]{OhWinter2014} for arbitrary moduli.

\section{Proofs}

The proofs are a relatively minor adaptation of the argument in \cite{MageeOhWinter2015}, which builds on the breakthrough in \cite{OhWinter2014} (the latter is itself based on key ideas in \cite{BourgainGamburdSarnak2011} combined with \cite{Dolgopyat1998, Naud2005, Stoyanov2011}).
%
%
The point of departure from the treatment in
\cite{MageeOhWinter2015} 
is
in the analysis of the measure
 $\mu_{s,x,\ga^{M}}$ in equation (135), culminating in Lemma 4.8, valid for arbitrary moduli $q$. We will follow this treatment,
henceforth  importing all the concepts and notation
from
 that paper.

Thus we are lead to study the measure $\mu$ on $G=\SL_{2}(q)$ given by
\be\label{eq:muIs}
\mu \ = \ \mu_{s,x,\ga^{M}}
\ \equiv
\
\sum_{\ga^{N} > \ga^{M}}
\exp(
[
\tau_{a}^{N}+ib\tau^{N}
]
(\ga^{N}x)
)
\gd_{c_{q}^{R}(\ga^{R}x)}
,
\ee
as in \cite[(135)]{MageeOhWinter2015}.
Here $x\in I$, $\ga^{M}$ is a fixed branch of $T^{-M}$, and $\ga^{N}=\ga^{M}\ga^{R}$. 
For ease of exposition, we first assume
  that we are treating the full shift as in \thmref{thm:2}, and that we can therefore view sums over branches $\ga^N$ as sums over globally (on $I$) defined branches of $T^{-N}$. 
Moreover, assume for simplicity that $\G(\mod q)=\SL_{2}(q)$. (Both of these assumptions are satisfied in the Zaremba setting of \cite{BourgainKontorovich2014}.)

Our goal in this paper is to prove the following
\begin{thm}\label{thm:main}
For $|a-s_{0}|<a_{0}$ and $\vf\in E_{q}$ (as defined in \cite[\S4.1]{MageeOhWinter2015}), we have
\be\label{eq:main}
\|\mu * \vf\|_{2} \ \le \ C\
q^{-1/4}\,
B\,
\|\vf\|_{2},
\ee
where
$$
\|\mu\|_{1}\ <\ B.
$$
\end{thm}

This is the replacement of \cite[Lemma 4.5]{MageeOhWinter2015} (bypassing the property (MIX)), and the rest of the proof of \cite[Lemma 4.8]{MageeOhWinter2015} follows analogously.

\

To begin, we 
 pick some $o \in I$,
and
define the measure $\nu$ by:
\be\label{eq:nuDef}
\nu \ \equiv \
\exp(\tau_{a}^{M}(\ga^{M}o))\mu_{1},
\ee
where $\mu_{1}$ is the measure given by
\be\label{eq:mu1Def}
\mu_{1} \ \equiv \
\sum_{\ga^{R}}
\exp(\tau_{a}^{R}(\ga^{R}o))
\gd_{c_{q}^{R}(\ga^{R}o)}
.
\ee

\begin{lem}
We have
\be\label{eq:muBndNu}
|\mu| \ \le \ C\, \nu.
\ee
\end{lem}
\pf
Use the ``contraction property'' in
\cite[(145-146)]{MageeOhWinter2015} and
argue
 as in the proof of \cite[Lemma 4.4]{MageeOhWinter2015}.
\epf

We will now manipulate $\mu_{1}$. We assume that $R$ can be decomposed further  as 
\be\label{eq:RpL}
R \ = \ R'L,
\ee
with $L$ to be chosen later (a sufficiently large constant independent of $R'$ and $q$).

Now split $\ga^R$ as
\begin{equation}\label{eq:arsplitting}
\ga^R \  = \  \ga_{R'}^L \ga_{R' - 1}^L \ldots \ga^L_2 \ga^L_1
,
\end{equation}
where the $\ga_{k}^L$ are  branches of $T^{-L}$.
This splitting \eqref{eq:arsplitting} is uniquely determined by $\ga^R$. For each $k\ge2$, we also split 
$$
\ga^L_k \  = \  \ga^{L-1}_k \ga^1_k  ,
$$
where $\ga^1_k = g_{i_k}$ for some $i_k$.

Write out
\bea
\nonumber
\tau_a^R(\ga^R o)  &=& \sum_{ i = 0 }^{R-1} \tau_a( T^i \ga^R o ) 
\\
\nonumber
&=& \sum_ {i = 0 }^{R'-1} \sum_{ \ell = 0}^{L-1} \tau_a( T^{  i L+ \ell } \ga^R o) 
\\
\nonumber
&=&  \sum_ {i = 0 }^{R'-1} \sum_{ \ell = 0}^{L-1} \tau_a( T^{  i L+ \ell } \ga^L_{R' - i } \ga^L_{R' - i - 1} \ldots \ga^L_1   o ) 
\\
\label{eq:lastline}
&=& \sum_{ i = 0}^{R' - 1} \tau_a^L (  \ga^L_{R' - i } \ga^L_{R' - i - 1} \ldots \ga^L_1 (o )) . 
\eea
We now perform decoupling term by term in the above.   We will use the shorthand
$$
\ga^{ L j  }  \ \equiv \  \ga^L_{j } \ga^{L}_{j - 1} \ldots \ga^L_1 .
$$
For $j \ge2$,
we  compare  each term
 in \eqref{eq:lastline}
  of the form
$$
\tau_a^L ( \ga^{ L j }( o ) )
$$
to 
$$
 \tau_a^L (  \ga^L_{j } \ga^{L-1}_{j-1} o ) 
. 
$$
This gives
\bea
\nonumber
\tau_a^L ( \ga^{ L j }( o ) ) &=&  \tau_a^L (  \ga^L_{j } \ga^{L-1}_{j-1} o ) + O \bigg( \sup | [ \tau_a^L \circ \ga^L_{j } ]' | d( \ga^{L-1}_{j- 1} o,   \ga^{L-1}_{j-1}\ga^1_{j-1} \ldots \ga^L_1 o )\bigg) 
\\
\label{eq:constant}
&=&  \tau_a^L (  \ga^L_{j } \ga^{L-1}_{j-1} o) + O ( \g^{-(L-1) } ) ,
\eea
where we used \cite[(68)]{MageeOhWinter2015}, valid when $a$ is suitably close to $s_0$. 

We will also use the formula
\begin{equation}\label{eq:cqsplit}
\delta_{c_q^R(\ga^R o) }  =  \delta_{c_q^L(\ga^L
o) }  *  \delta_{c_q^L(\ga^{2L} o) } *  \delta_{c_q^L(\ga^{3L} o) } * \ldots *      \delta_{c_q^L(\ga^{R' L} o) } .
\end{equation}
Then 
combining 
\eqref{eq:lastline} and \eqref{eq:cqsplit},
we write
\bea
\nonumber
\mu_1 &=&  \sum_{ \ga^{L}_1 , \ga^{L-1}_2 , \ldots , \ga^{L-1}_{R'}    } \sum_{ \ga^1_2 ,\ldots \ga^1_{R'}  }\exp( \tau_a^R  (\ga^R o ) )) \delta_{ c_q^R(\ga^R o)} 
\\
\nonumber
&=&  \sum_{ \ga^{L}_1 , \ga^{L-1}_2 , \ldots , \ga^{L-1}_{R'}    }  \sum_{ \ga^1_2 ,\ldots \ga^1_{R'}  } \exp \left ( \sum_{ j = 1}^{R'} \tau_a^L ( \ga^{ j L }(o))  \right)  
\times
\\
\label{eq:innersum}
&&
\hskip1in
\delta_{c_q^L(\ga^L o) }  *  \delta_{c_q^L(\ga^{2L} o) } *  \delta_{c_q^L(\ga^{3L} o) } * \ldots *      \delta_{c_q^L(\ga^{R' L} o) }  
.
\eea

We now decouple, replacing
 each term of the form
$$
e^{ \tau_a^L ( \ga^{ jL }(o) )  } \ \mapsto \ e^{ \tau_a^L (\ga^L_j \ga^{L-1}_{j-1} o) } 
 \ \equiv \ \gb_{j}
$$
with $j \ge2$,
at a cost of a multiplicative factor of $\exp(  c \g^{-L } )$; here $c$ is proportional to the implied constant of \eqref{eq:constant}. When $j = 1$, no replacement is performed, and we set
 $\gb_{1}\equiv e^{ \tau_a^L (\ga^L_{1} o)  }$. 

Inserting this into \eqref{eq:innersum} 
gives
\bea
\label{eq:mu1Bnd2}
\mu_{1}
&\leq &
\sum_{ \ga^{L-1}_1 , \ga^{L-1}_2 , \ldots , \ga^{L-1}_{R'}    }  
\sum_{\ga_{1}^{1}}
\gb_{1}
\delta_{c_q^L(\ga^L o) }  *
\\
\nonumber
&&
\exp( c \g^{-L} )^{R'-1} 
\left(
\sum_{ \ga^1_2 ,\ldots \ga^1_{R'}  } \prod_{ j =2 }^{R'} 
\gb_{j}\
\delta_{c_q^L(\ga^{2L} o) } *  \delta_{c_q^L(\ga^{3L} o) } * \ldots *      \delta_{c_q^L(\ga^{R' L} o) }  
\right).
\eea
Note that, although $\gb_{j}$ depends on all of the indices in $\ga_{j}^{L}\ga_{j-1}^{L-1}$,  
because
$\ga_{j}^{L-1}$ and $\ga_{j-1}^{L-1}$ are fixed in the outermost sum, we 
treat
$\gb_{j}$ as a function of  $\ga_{j}^{1}$. 

We claim that each term  $c_q^L(\ga^{jL } o)$ also only depends on one $\ga_j^1$. This is because 
we have $\ga^{jL} = g_{k_1} \ldots g_{k_L} \ga^{(j-1)L}$ for some choice of $g_{k_m}$, and
hence
 for whatever $o$ is chosen, we have
$$
c_q^L ( \ga^{jL } o)\  = \ c_q( g_{k_L}  \ga^{(j-1)L}o) c_q( g_{k_{L-1}} g_{k_L}  \ga^{(j-1)L} o ) \ldots  c_q( g_{k_1} \ldots g_{k_L}  \ga^{(j-1)L}o)
,
$$
see  \cite[(69)]{MageeOhWinter2015}.
Since $g_{k_m}$ maps 
$I$ into $I_{k_m}$, we have
$$
c_q( g_{ k_{m} } o' ) \ = \ g_{k_{m} } \bmod q
$$
for any $o' \in I$. Thus
\be\label{eq:cqL}
c_q^L( \ga^{jL} o )\  =\  g_{k_L} \ldots g_{k_1} \bmod q.
\ee
Here 
\be\label{eq:gjL}
g_{k_L} \ = \  \ga^1_j.
\ee

This means we may distribute the convolution and product over the sum,
writing \eqref{eq:mu1Bnd2} as
\bea
\nonumber
\mu_{1}
&\leq& 
\exp( c \g^{-L} )^{R'-1}
\sum_{ \ga^{L-1}_1 , \ga^{L-1}_2 , \ldots , \ga^{L-1}_{R'}    }  
\left(
\sum_{\ga_{1}^{1}}
\gb_{1}\
\delta_{c_q^L(\ga^L o) }  
\right)*
 \left( 
\sum_{ \ga_2^1 } 
\gb_{2}\
 \delta_{c_q^L(\ga^{2L} o) } 
\right) * \ldots
\\
\label{eq:convolve}
&&
\hskip2in
\ldots
 * \left( \sum_{\ga_{R'}^1 }
  \gb_{R'}\
  \delta_{c_q^L(\ga^{R'L} o) } \right).
\eea
We give each convolved term in \eqref{eq:convolve} a name, defining, for each $j\ge1$, the measure
\be\label{eq:etaDef}
\eta_{j} \ = \ \eta_{j}^{(\ga_{j}^{L-1},\ga_{j-1}^{L-1})}
\
\equiv
\
\sum_{ \ga_j^1 }
 \gb_{j }\
   \delta_{c_q^L(\ga^{jL} o) } 
.
\ee
We have thus  proved the following
\begin{prop}
We have
\be
\label{eq:mu1Prop}
\mu_{1}
\ \le \
\exp( c \g^{-L} )^{R'-1} 
  \sum_{ \ga^{L-1}_1 , \ga^{L-1}_2 , \ldots , \ga^{L-1}_{R'}    } 
 \eta_{1}
 * 
\eta_{2}
* \ldots
 *
\eta_{R'}
  .
\ee 
\end{prop}

Next we observe that each of the measures $\eta_{j}$ is nearly flat, in that their coefficients
in \eqref{eq:etaDef}
 differ by constants:
\begin{lem}
For each $j\ge1$ and any $\ga_{j}^{1}$ and ${\ga_{j}^{1}}'$, we have
\be\label{eq:lem1}
{\gb_{j}'\over \gb_{j}}
\ \le\ 
\exp(c\g^{-L+1})
.
\ee
\end{lem}
\pf
The first $L-1$ terms of $\gb_{j}$ and $\gb_{j}'$ agree, so we 
again use the ``contraction property'' 
\cite[(145-146)]{MageeOhWinter2015}.
\epf

Since the measures $\eta_{j}$ are nearly flat, we may now apply the expansion result in \cite{BourgainVarju2012}.

\begin{thm}\label{thm:flatExpand}
Assume $L$ is sufficiently large (depending only on $\G$). Then
for $\varphi\in L^{2}_{0}(G)$, we have
\be\label{eq:prop1}
\| \eta_{j}*\varphi \| _{2}
\ \le \
(1-C_{1})\,
\|\eta_{j}\|_{1}\,
\|\varphi\|_{2},
\ee
Here $C_{1}>0$ depends on $\G$  but not on $q$.
\end{thm}

To prove this theorem, we 
need the following simple
\begin{lem}\label{lem:Expand}
Let $\pi$ be a
unitary
 $G$-representation on a Hilbert space $\cH$, and assume 
that the operator $A$ acts on $\cH$ via
$$
A\vf \ =  \ \sum_{j\in J}\pi(h_{j})\vf,
$$
for some $h_{j}\in G$ and indexing set $J$. Assume that $A$ has the ``spectral gap'' property: there is some $C_{0}>0$ so that
\be\label{eq:Aacts}
\<A\vf,\vf\> \ \le \ (1-C_{0})\, |J| \, \|\vf\|^{2}.
\ee
For some positive coefficients $\gk_{j}>0$, let $\widetilde A$ act on $\cH$ as
$$
\widetilde A\vf \ = \ \sum_{j\in J}\gk_{j}\ \pi(h_{j})\vf,
$$
and assume that the $L^{\infty}$ norm of the coefficients is controlled by the $L^{1}$ norm, in the sense that for some $K\ge1$,
\be\label{eq:maxK}
\max \gk_{j} \ \le \ K \,\bar\gk,
\ee
where
$$
\bar\gk\ := \ \frac1{|J|}\sum_{j}\gk_{j}
$$
is the coefficient average.
Then $\widetilde A$ has the following ``spectral gap'':
\be\label{eq:AtilActs}
\<\widetilde A\vf,\vf\> \ \le \  \bar\gk\ (1-C_{0}+\sqrt{K-1})\, |J| \, \|\vf\|^{2}.
\ee
\end{lem}
\pf
This is an exercise in Cauchy-Schwarz.
\epf

With this lemma, it is a simple matter to give a

\pf[Proof of \thmref{thm:flatExpand}]\

We will apply \lemref{lem:Expand} with $\cH=L^{2}_{0}(G)$ and $\pi$ the right-regular representation. 
Recalling \eqref{eq:etaDef}, 
we can write
$$
\| \eta_{j}*\varphi \|^{2} _{2}\  = \ \<\widetilde A \vf,\vf\>
,
$$
where $\widetilde A
$ acts by convolution with the measure
$$
\sum_{a_{j}^{1},{a_{j}^{1}}'}
\gb_{j}\,\gb_{j}'\,
\delta_{c_q^L(\ga^{jL} o) c_q^L((\ga^{jL})' o)^{-1} } .
$$
Using the notation of \eqref{eq:cqL} and \eqref{eq:gjL}, note that
$$
c_q^L(\ga^{jL} o) c_q^L((\ga^{jL})' o)^{-1} 
\ = \
\ga_{j}^{1}\cdot g_{k_{L-1}} \ldots g_{k_1}
({\ga_{j}^{1}}'\cdot g_{k_{L-1}} \ldots g_{k_1})^{-1}
\ = \
\ga_{j}^{1}
({\ga_{j}^{1}}')^{-1}.
$$
The indexing set $J$ of \lemref{lem:Expand} then runs over pairs $\ga_{j}^{1},{\ga_{j}^{1}}'$,
 the coefficients $\gk_{j}$ are the products $\gb_{j}\gb_{j}'$,
 and the elements $h_{j}$ are 
 $\ga_{j}^{1}
({\ga_{j}^{1}}')^{-1}$.

That the operator $A$ (without coefficients) has a spectral gap \eqref{eq:Aacts} is precisely the statement proved in \cite{BourgainVarju2012}, with $C_{0}$ independent of $q$.\footnote{%
Here we need the products  $\ga_{j}^{1}
({\ga_{j}^{1}}')^{-1}$ to generate group with Zariski closure $\SL_{2}$. In the Zaremba case,
it is important
 that each $\ga_{j}^{1}$ is a product of two generators $\mattwos011a\mattwos011b$. Otherwise, e.g., the products $\mattwos011a\mattwos011b^{-1}=\mattwos10{a-b}1$ could all be lower-triangular.} 
The bound \eqref{eq:maxK} follows from \eqref{eq:lem1} with
$$
K \ = \
\exp(2c\g^{-L+1}).
$$ 
Note also that
$$
|J|\bar \gk \ = \
\left(\sum_{\ga_{j}^{1}}\gb_{j}\right)^{2}
\ = \
\|\eta_{j}\|_{1}^{2}
.
$$
Choosing $L$ sufficiently large (depending only on $\G$), one can make $K$ sufficiently close to $1$ so that \eqref{eq:AtilActs} gives \eqref{eq:prop1}, as claimed.
\epf

\begin{cor}
Assume that $L$ is sufficiently large (depending only on $\G$). Then there is some $C_{2}>0$ also depending only on $\G$ so that, for any $\vf\in L^{2}_{0}(G)$, we have
\be\label{eq:mu1Bnd}
\|\mu_{1}*\vf\|_{2} \ \le \ 
(1-C_{2})^{R} \ \|\mu_{1}\|_{1} \ \|\vf\|_{2}.
\ee
\end{cor}
\pf
Beginning with \eqref{eq:mu1Prop}, apply \eqref{eq:prop1} $R'$ times to get
$$
\|\mu_{1}*\vf\|_{2} \ \le \ 
\exp(c\g^{-L})^{R'-1}
\sum_{\ga_{1}^{L-1},\dots,\ga_{R'}^{L-1}}
(1-C_{1})^{R'}
\prod_{j=1}^{R'}
\|\eta_{j}\|_{1}\|\vf\|_{2}.
$$
Applying contraction yet again gives
$$
\sum_{\ga_{1}^{L-1},\dots,\ga_{R'}^{L-1}}
\prod_{j=1}^{R'}
\|\eta_{j}\|_{1}
\ \le \
\exp(c\g^{-L})^{R'-1}
\|\mu_{1}\|_{1},
$$
whence \eqref{eq:mu1Bnd} follows on taking $L$ large enough and recalling  \eqref{eq:RpL}.
\epf

Returning to the measure $\nu$ in \eqref{eq:nuDef},
we have from \eqref{eq:mu1Bnd} that
\be\label{eq:nuBndR}
\|\nu*\vf\|_{2} \ \le \ 
(1-C_{2})^{R} \ \|\nu\|_{1} \ \|\vf\|_{2}.
\ee
To conclude \thmref{thm:main}, we need the following
\begin{lem}
Let $\mu$ be a complex distribution on $G=\SL_{2}(q)$ and assume that $|\mu|\le C\nu$. Let $E_{q}\subset L_{0}^{2}(G)$ be the subspace defined in \cite[\S4.1]{MageeOhWinter2015}, and let $A:E_{q}\to E_{q}$ be the operator acting by convolution with $\mu$. Then
\be\label{eq:1p5}
\| A\| \ \le\ C' \left[{|G| \ \|\widetilde\nu * \nu\|^{2}_{2}\over q}\right]^{1/4}.
\ee
Here $\widetilde \mu(g)=\overline{\mu(g^{-1})}$.
\end{lem}
\pf
Note that the operator $A^{*}A$ is self-adjoint, positive, and acts by convolution with $\widetilde\mu*\mu$. Let $\gl$ be an eigenvalue of $A^{*}A$. Since $A$ acts on $E_{q}$, Frobenius gives that $\gl$ has multiplicity $\mult(\gl)$ at least $Cq$. We then have that
\beann
\gl^{2}\ \mult(\gl) 
&\le& 
\tr[(A^{*}A)^{2}]
\ = \
\sum_{g\in G}
\<(A^{*}A)^{2}\gd_{g},\gd_{g}\>
\ = \
\sum_{g\in G}
\|\widetilde \mu * \mu *\gd_{g}\|_{2}^{2}
\\
& =&
|G|\
\|\widetilde \mu * \mu \|_{2}^{2}
\ \le \
C^{4}\ 
|G|\
\|\widetilde \nu * \nu \|_{2}^{2}.
\eeann
The claim follows, as $\|A\|=\max_{\gl}\gl^{1/2}$.
\epf

We apply the lemma to $\mu$ in \eqref{eq:muIs}
using
\eqref{eq:muBndNu}, giving
\be\label{eq:muConvBnd}
\|\mu * \vf\|_{2}\ \le\ C\,q^{1/2}
 \|\widetilde\nu * \nu\|^{1/2}_{2}
 .
\ee
It remains to estimate the $\nu$ convolution.
\begin{prop}
Choosing $R$ to be of size $C\log q$ for suitable $C$, we have that
\be\label{eq:nuBnd}
 \|\widetilde\nu * \nu\|_{2}
 \
 \le
 \
 2{\|\nu\|_{1}^{2}\over |G|^{1/2}}.
\ee
\end{prop}
\pf
Let
$$
\psi \ \equiv \ \gd_{e} - \frac1{|G|}\bo_{G} \ \in \ L_{0}^{2}(G),
$$
and note that 
$
\|\psi\|_{2}< 1.
$ 
Then
\beann
 \|\widetilde\nu * \nu\|_{2}
 &=&
  \|\widetilde\nu * \nu*\gd_{e}\|_{2}
\ \le\
\|\widetilde\nu * \nu*\left(\frac1{|G|}\bo_{G}\right)\|_{2}
+
\|\widetilde\nu * \nu*\psi\|_{2}
\\
& \le&
{\|\nu\|_{1}^{2}\over |G|^{1/2}}
+
\|\nu\|_{1}
\| \nu*\psi\|_{2}
,
\eeann
where we used the triangle inequality and Cauchy-Schwarz. 
Since $\psi\in L_{0}^{2}(G)$, we apply \eqref{eq:nuBndR}, giving
$$
\| \nu*\psi\|_{2}
<
(1-C_{2})^{R} \ \|\nu\|_{1}
<
{\|\nu\|_{1}\over |G|^{1/2}}
$$
by a suitable choice of $R=C\log q$. The claim follows immediately.
\epf

Finally, we give a
\pf[Proof of \thmref{thm:main}]
Insert \eqref{eq:nuBnd} into \eqref{eq:muConvBnd} and use \eqref{eq:muBndNu} and $|G|>Cq^{3}$. Clearly \eqref{eq:main} holds with $B=C\|\nu\|_{1}$. 
\epf


\subsection{Modifications for Subshifts}\label{schottky}

We sketch here the modifications needed to handle the case $\G$ is a Schottky group as in \thmref{thm:1}. Then $I=\cup_k I_k$, where to each  $I_k$ is assigned some $g_k\in\SL_2(\Z)$ such that $T|_{I_k}=g_k^{-1}$ and $c_0|_{I_k}\equiv g_k$. The shift is restricted to exclude any letter $g_k$ being followed by $g_k^{-1}$. Note that while in \cite{MageeOhWinter2015} it is stated that the values $c_0(I)$ should freely generate a semigroup, the arguments also apply equally to the Schottky case. 

In the decomposition \eqref{eq:innersum}, each sum on $\ga_j^1$ needs to be restricted to be admissible, once $\ga_{j-1}^{L-1}$ and $\ga_{j}^{L-1}$ are chosen (and each itself is an admissible sequence). 
The base points $o\in I$ need to be chosen in the appropriate domains of branches of $T^{-L}$, etc.; we only ever use the contraction principle, so  these choices have no effect.

The following  issue arises when  $\G$ is generated by two elements, $g$ and $h$, say. Suppose $\ga_{j-1}^{L-1}$ ends in $g$ while $\ga_j^{L-1}$ starts with $g^{-1}$. Then in the $\ga_j^1$ sum, only $h$ and $h^{-1}$ are admissible, and this does not generate a Zariski dense group for the operator $A$ in the proof of \thmref{thm:flatExpand}. To fix this issue, one instead decomposes each block $\ga_j^L$ as $\ga_j^{L-2}\ga_j^2$, that is, isolating two indices instead of one. 
With this adjustment, even if
$\ga_j^{L-2}$ ends
 in $g$ and $\ga_{j-1}^{L-2}$ starts
 in $g^{-1}$, 
the admissible $\ga_j^2$ sum runs over the elements
$
gh,
gh^{-1},
hg^{-1},
hh,
h^{-1}g^{-1},
h^{-1}h^{-1}
$.
It is then easy to see that the operator $A$ in the proof of \thmref{thm:flatExpand} generates a Zariski dense group
(if 
$\G$ has more than two generators, this
 is 
 clear). Now, this group and its generator set (and hence also its expansion constant $C_0$ in \eqref{eq:Aacts}) depend on $\ga_j^{L-2}$ and $\ga_{j-1}^{L-2}$ (or rather just their starting/ending letters). But as $\G$ is finitely generated, only a finite number of groups/generators arise in this way, and
we simply
take $C_0$ to be the worst 
one%
. With these modifications, the proof goes through as before.

\bibliographystyle{alpha}

\bibliography{../../AKbibliog}

\begin{thebibliography}{MOW15}

\bibitem[BGS11]{BourgainGamburdSarnak2011}
J.~Bourgain, A.~Gamburd, and P.~Sarnak.
\newblock Generalization of {S}elberg's 3/16th theorem and affine sieve.
\newblock {\em Acta Math}, 207:255--290, 2011.

\bibitem[BK14]{BourgainKontorovich2014}
J.~Bourgain and A.~Kontorovich.
\newblock On {Z}aremba's conjecture.
\newblock {\em Annals Math.}, 180(1):137--196, 2014.

\bibitem[BV12]{BourgainVarju2012}
Jean Bourgain and P{\'e}ter~P. Varj{\'u}.
\newblock Expansion in {$SL_d({\bf Z}/q{\bf Z}),\,q$} arbitrary.
\newblock {\em Invent. Math.}, 188(1):151--173, 2012.

\bibitem[Dol98]{Dolgopyat1998}
Dmitry Dolgopyat.
\newblock On decay of correlations in {A}nosov flows.
\newblock {\em Ann. of Math. (2)}, 147(2):357--390, 1998.

\bibitem[MOW15]{MageeOhWinter2015}
Michael Magee, Hee Oh, and Dale Winter.
\newblock Expanding maps and continued fractions, 2015.
\newblock Preprint, {\tt arXiv:1412.4284v2}.

\bibitem[Nau05]{Naud2005}
Fr{\'e}d{\'e}ric Naud.
\newblock Expanding maps on {C}antor sets and analytic continuation of zeta
  functions.
\newblock {\em Ann. Sci. \'Ecole Norm. Sup. (4)}, 38(1):116--153, 2005.

\bibitem[OW14]{OhWinter2014}
Hee Oh and Dale Winter.
\newblock Uniform exponential mixing and resonance free regions for convex
  cocompact congruence subgroups of {S}{L}(2,{Z}), 2014.
\newblock Preprint, {\tt arXiv:1410.4401v2}.

\bibitem[Sto11]{Stoyanov2011}
Luchezar Stoyanov.
\newblock Spectra of {R}uelle transfer operators for axiom {A} flows.
\newblock {\em Nonlinearity}, 24(4):1089--1120, 2011.

\end{thebibliography}

\end{document}